\documentclass[11pt]{article}

\usepackage{float}
\usepackage{fullpage}
\usepackage{amsmath,amssymb,amsthm}
\usepackage{marginnote} 
\usepackage[round]{natbib}
\setcitestyle{numbers,square}
\usepackage{graphicx}
\usepackage{caption}
\usepackage{authblk}

\numberwithin{equation}{section}
\newtheorem{theorem}{Theorem}

\newtheorem{proposition}{Proposition}[section]

\newtheorem{remark}{Remark}[section]

\usepackage{commath}
\usepackage[normalem]{ulem}
\usepackage{mathabx}
\usepackage{mathrsfs}

\usepackage[usenames]{color}
\definecolor{plum}{rgb}{.4,0,.4}
\definecolor{BrickRed}{rgb}{0.6,0,0}
\definecolor{light-gray}{gray}{0.5}
\usepackage[plainpages=false,pdfpagelabels,colorlinks=true,linkcolor=BrickRed,citecolor=plum]{hyperref}

\def\ddefloop#1{\ifx\ddefloop#1\else\ddef{#1}\expandafter\ddefloop\fi}

\def\ddef#1{\expandafter\def\csname c#1\endcsname{\ensuremath{\mathcal{#1}}}}
\ddefloop ABCDEFGHIJKLMNOPQRSTUVWXYZ\ddefloop

\def\ddef#1{\expandafter\def\csname s#1\endcsname{\ensuremath{\mathsf{#1}}}}
\ddefloop ABCDEFGHIJKLMNOPQRSTUVWXYZ\ddefloop

\def\E{\mathop{\mathbf{E}}}

\begin{document}

\title{Estimating Certain Integral Probability Metric (IPM) Is as Hard as Estimating under the IPM}

\author[1]{Tengyuan Liang\thanks{Liang gratefully acknowledges support from the George C. Tiao Fellowship. The paper was previously titled ``On the minimax optimality of estimating Wasserstein metric'' \citep{liang2019minimax}. The previous version is no longer intended for publication.}}

\affil[1]{University of Chicago}
\date{}
\maketitle

\begin{abstract}
	We study the minimax optimal rates for estimating a range of Integral Probability Metrics (IPMs) between two unknown probability measures, based on $n$ independent samples from them. Curiously, we show that estimating the IPM itself between probability measures, is not significantly easier than estimating the probability measures under the IPM. We prove that the minimax optimal rates for these two problems are multiplicatively equivalent, up to a $\log \log (n)/\log (n)$ factor.
\end{abstract}

\section{Introduction}

In this note, we study the minimax optimal rates for estimating the Integral Probability Metrics (IPMs) between probability measures based on samples. IPMs are widely used in both statistics and machine learning, with applications in nonparametric two-sample tests \citep{sriperumbudur2012empirical,gretton2012kernel}, inferring the transportation cost (the Wasserstein-1 metric) from one set of samples to another \citep{sommerfeld2018inference,peyre2019computational}, and with more recent appearances in rigorous investigations on the generative adversarial networks (GANs) \citep{arjovsky2017wasserstein, liu2017approximation, liang2018well,singh2018nonparametric,uppal2019nonparametric,binkowski2018demystifying}.

Let $\mu, \nu$ be two probability measures supported on $\Omega=[0, 1]^d$, and $d_{\cF}(\mu, \nu)$ denote a certain IPM between them induced by a set of functions $\cF$, defined as
\begin{align}
	d_{\cF}(\mu, \nu) := \sup_{f \in \cF} | \int_{\Omega} f d\mu - \int_{\Omega} fd\nu |\enspace.
\end{align} 
Consider that $X_1, \ldots X_m$ are i.i.d samples from $\mu$, and $Y_1, \ldots, Y_n$ i.i.d from $\nu$. We study the minimax optimal rate for estimating $d_{\cF}(\mu, \nu)$ based on $\{X_i\}_{i=1}^m, \{Y_j\}_{j=1}^n$, for some class of probability measures $\cG$ of interest
    \begin{align}
        \inf_{\widetilde T_{m,n}} \sup_{\mu, \nu \in \cG} \E | \widetilde T_{m,n} -  d_{\cF}(\mu, \nu) | \enspace.
    \end{align}

It turns out that using the empirical measure $\widehat \nu_n$ to estimate is a bad idea when $\cF$ is complex enough, regardless of how simple $\cG$ is. To see this, let's consider a simple case with $\cF = {\rm Lip}(1)$. In such a case, $d_{\cF}$ reduces to the Wasserstein-1 metric $W$ \eqref{eq:wass}. Due to a result by \citet{dudley1969speed}, even for infinitely smooth $\cG = \{{\rm Unif}(\Omega) \}$ and $d\geq 2$,
\begin{align}
     \sup_{ \nu \in \cG} |W(\mu, \widehat\nu_n) -  W(\mu, \nu)| \asymp  n^{-\frac{1}{d}} \enspace.
\end{align}
A natural question arises: can one obtain faster rates, for estimating the IPM with other estimators $\widetilde T_{m,n}$ leveraging certain regularity of $\cG$ such as smoothness?

A related yet different problem studied in the current literature is estimating a probability measure under certain IPMs \citep{sriperumbudur2012empirical, weed2017sharp, liang2018well, singh2018nonparametric, weed2019estimation}, in the following sense
    \begin{align}
        \inf_{\widetilde\nu_n} \sup_{\nu \in \cG} \E d_{\cF}(\widetilde\nu_n , \nu) \enspace.
    \end{align}
The two problems are closely related: ``estimating the metric itself'' is usually an \textbf{easier} problem than ``estimating under the metric.'' In fact, the solution of the latter problem $\widetilde \mu_m, \widetilde \nu_n$ naturally induces a plug-in answer to the former, since
\begin{align*}
    \E |d_{\cF}(\widetilde \mu_m, \widetilde \nu_n) - d_{\cF}(\mu, \nu)| \leq \E d_{\cF}(\widetilde\mu_m , \mu) + \E d_{\cF}(\widetilde\nu_n , \nu) \enspace.
\end{align*}
However, it is unclear whether such a plug-in estimator is optimal. In fact, it is well-known that estimating specific functional of density $F(\nu)$ is usually strictly easier than estimating the density $\nu$ itself. For example, in estimating quadratic functionals of a smooth density vs. estimating under the quadratic functionals, the plug-in approach is strictly sub-optimal, where the rates can be much-improved \citep{bickel1988estimating,donoho1990minimax,fan1991estimation}. In recent practical applications such as GANs, one is curious to understand if evaluating and inferring how well we do in terms of learning the probability measure, could be simpler than learning the measure itself \citep{lucic2018gans, liang2018interaction}.

In this paper, however, we prove that ``estimating the IPMs,'' is \textbf{not significantly easier} than ``estimating under the IPMs,'' for a wide range of measures and metrics. Specifically, the plug-in approach is minimax optimal up to a $\log \log (n)/\log (n)$ factor
\begin{align*}
    \frac{\log \log (n \wedge m)}{\log (n \wedge m)} \cdot (n \wedge m)^{-\frac{\beta+\gamma}{2\beta+d}}  \precsim \inf_{\widetilde T_{m,n}} \sup_{\mu, \nu \in \cG_\beta} &\E | \widetilde T_{m,n} -  d_{\cF_\gamma}(\mu, \nu) |  \\
     \leq  \inf_{\widetilde\mu_m, \widetilde\nu_n} \sup_{\mu, \nu \in \cG_\beta} &\E |d_{\cF_\gamma}(\widetilde \mu_m, \widetilde \nu_n) - d_{\cF_\gamma}(\mu, \nu)|  \precsim (n \wedge m)^{-\frac{\beta+\gamma}{2\beta+d}} .
\end{align*}
Here $\cG_\beta$ contains probability measures with densities in the H\"{o}lder space with smoothness $\beta \in \mathbb{R}_{\geq 0}$, and the IPMs are induced by $\cF_\gamma$, the H\"{o}lder space with smoothness $\gamma \in \mathbb{R}_{\geq 0}$, with $\gamma < d/2$. Note that when $\gamma \geq d/2$, the parametric rate $n^{-1/2}$ is attainable.
The result informs us that (1) seeking for other forms of estimators for $d_{\cF_\gamma}(\mu,\nu)$ would only improve the rates logarithmically, and (2) estimating the IPM between two measures is fundamentally just as hard as estimating the measure under the IPM.

\subsection{Preliminaries}

We introduce the notations used in the paper.
For a function $f:\mathbb{R}^d \rightarrow \mathbb{R}$ and $p\geq 1$, $\| f \|_{L_p}$ denotes the $L_p$ norm w.r.t. the Lebesgue measure. For a finite dimensional vector $\theta$ and $q\geq 1$, $\| \theta \|_q$ is the vector $L_q$ norm, and $\| \theta \|$ is the $L_2$ norm. For an integer $K$, $[K]:=\{0, 1,\ldots, K-1\}$.

Let $\sC^{\beta}:=\sC^{\lfloor \beta \rfloor, \beta-\lfloor \beta\rfloor}$ to be H\"{o}lder space with smoothness $\beta >0$. 
\begin{align}
    \sC^{\beta} :=\left\{ f: \Omega \rightarrow \mathbb{R}: \max_{|\alpha| \leq  \lfloor \beta \rfloor}  \sup_{x \in \Omega} |D^{\alpha} f|  +  \max_{|\alpha| =  \lfloor \beta \rfloor}\sup_{x \neq y \in \Omega}\frac{| D^{\alpha} f(x) - D^{\alpha} f(y) |}{\|x - y \|^{\beta-\lfloor \beta\rfloor}} <\infty \right\}
\end{align}
where $\alpha=[\alpha_1,\ldots,\alpha_d] \in \mathbb{N}^d$ ranges over multi-indices, and $|\alpha| := \sum_{i=1}^d \alpha_i$.
We only consider the bounded case with $\Omega = [0, 1]^d$.

Define the forward difference operator, given a function $f: \mathbb{R}^d \rightarrow \mathbb{R}$, for every $h \in \mathbb{R}^d$
\begin{align*}
    \Delta_h f(x) &:= f(x+h) - f(x), \\
    \Delta_h^m f(x) &:= \Delta_h ( \Delta_h^{m-1} f(x) ), ~m\geq 2.
\end{align*}
Let $1\leq p, q \leq \infty$ and $\beta>0$, the Besov-Lipschitz space semi-norm $\| \cdot \|^{'}_{\sB^{\beta, p}_{q}}$ is defined in the following way \cite[Chapter 17, Proposition 17.21]{leoni2017first}.
\begin{align}
    \| f \|^{'}_{\sB^{\beta, p}_{q}} := \left( \sum_{j=0}^{\infty}\left((2^j)^{\beta} \sup_{\| h\| \leq 1/2^j} \| \Delta_h^{\lfloor \beta \rfloor + 1 } f \|_{L_p} \right)^q \right)^{1/q} \enspace.
\end{align}

Wavelets are used to provide an equivalent characterization of the Besov spaces, if the basis $\{ h_{jk}, j \in \mathbb{N}, 0\leq k <2^{jd} \}$ satisfies certain regularity conditions \cite[Chapter 9, Theorem 9.2]{hardle2012wavelets}. For function $f(x) = \sum_{j=0}^{\infty} \sum_{k<2^{dj}} \theta_{jk} h_{jk}(x)$, define the Besov space norm in terms of the wavelet coefficients \citep{weed2019estimation}
\begin{align}
    \| f \|_{\sB^{\beta, p}_{q}} := \left( \sum_{j=0}^\infty \left((2^j)^{\beta} (2^{jd})^{\frac{1}{2}-\frac{1}{p}}  \| \theta_{j \cdot} \|_{p} \right)^q \right)^{1/q} \enspace .
\end{align}
In this paper, we assume such regularity conditions throughout so that the Besov space norms $\| \cdot \|^{'}_{\sB^{\beta, p}_{q}}$ and $\| \cdot \|_{\sB^{\beta, p}_{q}}$ are equivalent. We refer the readers to \cite[Appendix C]{weed2019estimation} and \cite[Chapter 2.12]{cohen2003numerical} for details on the regularity conditions that we assume on the wavelets.

Besov spaces subsume H\"{o}lder spaces as special cases ($p=q=\infty$), in the following sense
\citep{tribel1980theory,donoho1996density, hardle2012wavelets}:
        under regularity conditions, the following equivalence holds between the Besov space and H\"{o}lder space
        $$\sB^{\beta, \infty}_{\infty} = \sC^{\beta}, \text{for}~\beta \notin \mathbb{N}.$$
        In particular, when $\beta=1$, $\sB^{\beta, \infty}_{\infty}$ is called the Zygmund space, which contains the Lipschitz space
        $\sB^{1,\infty}_{\infty} \supseteq {\rm Lip} \supseteq \sB^{1,\infty}_{1}$.

Now we are ready to formally state the parameter spaces, and the IPMs to study.

\paragraph{Parameter Spaces.}
For some $M>0$, the class of probability measures of interest is 
\begin{align}
    \label{eq:measure-class}
    \cG_\beta:=\left\{ \mu ~:~ \int_\Omega d \mu = 1, \mu \geq 0, \frac{d\mu}{dx} \in \sB^{\beta,\infty}_{\infty}(M)   \right\} \enspace,
\end{align}
with
\begin{align*}
    \sB^{\beta, \infty}_{\infty}(M) := \left\{ f ~:~  \| f \|_{\sB^{\beta, \infty}_{\infty}} \leq M \right\} \enspace.
\end{align*}
Again, for non-integer $\beta$, we are considering densities that are H\"{o}lder smooth.

\paragraph{Integral Probability Metric.}
The class of IPMs considered is induced by the Besov space, for some $\gamma > 0$
\begin{align}
    \cF_\gamma &:= \sB^{\gamma, \infty}_{\infty}(1) \enspace, \nonumber\\
    \label{eq:ipm-class}
    d_{\cF_\gamma}(\mu, \nu) &= \sup_{f \in \cF_\gamma} | \int f d \mu - \int f d\nu | \enspace. 
\end{align}
As a special case for the IPMs, the Wasserstein-$1$ metric (for measures supported on bounded $\Omega$) is
\begin{align}
    \label{eq:wass}
    W(\mu, \nu) &:= \sup_{f \in {\rm Lip}(1)} | \int f d\mu - \int f d\nu| \enspace.
\end{align}

\section{Optimal Rates for Estimating IPMs}

\begin{theorem}[Minimax Rate]
    Consider the domain $\Omega = [0, 1]^d$. Given $m$ i.i.d. samples $X_1, \ldots, X_m$ from $\mu$, and $n$ i.i.d. samples $Y_1, \ldots, Y_n$ from $\nu$. Then the minimax optimal rate for estimating $d_{\cF_\gamma}(\mu,\nu)$ satisfies
    \begin{align}
        \frac{\log \log (n \wedge m)}{\log (n \wedge m)} \cdot (n \wedge m)^{-\frac{\beta+\gamma}{2\beta+d}}  \precsim \inf_{\widetilde T_{m,n}} \sup_{\mu, \nu \in \cG_\beta} &\E | \widetilde T_{m,n} -  d_{\cF_\gamma}(\mu, \nu) | \precsim (n \wedge m)^{-\frac{\beta+\gamma}{2\beta+d}} \enspace.
    \end{align}
	Here $\mu, \nu$ lie in $\cG_\beta, \beta \in \mathbb{R}_{\geq 0}$ as in \eqref{eq:measure-class} whose densities are $\beta$-H\"{o}lder smooth. The function class $\cF_\gamma, \gamma \in \mathbb{R}_{\geq 0}$ for the metric is defined in \eqref{eq:ipm-class}, with $\gamma < d/2$.
\end{theorem}
\begin{remark}
\rm    
Here the $\beta$ quantifies the regularity of the measures, and $\gamma$ quantifies the regularity of the metrics. 

A few remarks are in order. First, we emphasize that the main technicality is in deriving the lower bound. We construct two composite/fuzzy hypotheses using delicate priors with matching $\log (n \wedge m)$ moments. However, the IPMs to estimate differs sufficiently under the null vs. the alternative. Then we calculate the Total Variation (TV) metric directly on the posterior of data samples defined by the composite hypothesis, using some telescoping techniques involving sum-products. The transparent technique could be of independent interest in handling TV-type calculations in proving lower bounds.
Second, as a direct corollary, the following extension holds true. Suppose $\mu \in \cG_{\beta_1}$ and $\nu \in \cG_{\beta_2}$, then define $\beta :=\beta_1 \wedge \beta_2$,
\begin{align*}
    \frac{\log \log (n \wedge m)}{\log (n \wedge m)} \cdot (n \wedge m)^{-\frac{\beta+\gamma}{2\beta+d}}  \precsim \inf_{\widetilde T_{m,n}} \sup_{\mu \in \cG_{\beta_1}, \nu \in \cG_{\beta_2}} &\E | \widetilde T_{m,n} -  d_{\cF_\gamma}(\mu, \nu) | \precsim (n \wedge m)^{-\frac{\beta+\gamma}{2\beta+d}} \enspace.
\end{align*}
Third, the $\gamma < d/2$ condition is effectively equivalent to that $\cF_\gamma$ is beyond the Donsker's class. This is the complex regime since within the Donsker's class $\gamma \geq d/2$, the parametric rate $n^{-1/2}$ is attainable \citep{gretton2012kernel,liang2018well, singh2018nonparametric}.
Finally, we would like to remark that in a concurrent work, \citet{niles2019estimation} obtains a similar lower bound (an improvement over our result with a $(\log n)^{-1/d}$ factor) for the special case $\beta = 0$ and $\gamma = 1$, using a distinct approach. We agree that closing the $\log n$ gap is an interesting question for future work.

\end{remark}

\subsection{Proof of the Lower Bound}

Without the loss of generality, consider the case when $m \geq n$. The lower bound construction is divided into six logical steps, for better organization. We make use of multi-resolution analysis in the construction.

\paragraph{Step 1: reduction to Besov space semi-norm.}
For any $p \geq 1$, define $p_\star \geq 1$ such that $1/p_\star + 1/p = 1$. For simplicity, define Radon-Nikodym derivative of measure $\mu$ w.r.t. the Lebesgue measure as $\rho_\mu(x) := d\mu/dx$.
Define explicitly the wavelet coefficients $f_{jk} := \langle f, h_{jk} \rangle$, and $u_{jk} := \langle d \mu/dx, h_{jk}\rangle, v_{jk} := \langle d \nu/dx, h_{jk}\rangle$. Under such notations, the integral probability metric reduces to the following
\begin{align*}
	d_{\sB^{\gamma,p}_q}(\mu, \nu) &:= \sup_{f \in \sB^{\gamma,p}_{q}(1)} |  \int f d\mu  - \int f d\nu |\\
	& = \sup_{f \in \sB^{\gamma,p}_q(1)} | \sum_{j\geq 0} \sum_{k=0}^{2^{dj}-1} f_{jk} (u_{jk} - v_{jk}) | \\
	& = \sup_{f \in \sB^{\gamma,p}_q(1)}  | \sum_{j\geq 0} \| f_{j\cdot} \|_{p} \| u_{j\cdot} - v_{j\cdot} \|_{p_{\star}} | \\
	& = \sup_{f \in \sB^{\gamma,p}_q(1)}  | \sum_{j\geq 0} (2^{dj})^{\frac{\gamma}{d}+\frac{1}{2}-\frac{1}{p}} \| f_{j\cdot} \|_{p} \cdot (2^{-dj})^{\frac{\gamma}{d}+\frac{1}{2}-\frac{1}{p}} \| u_{j\cdot} - v_{j\cdot} \|_{p_{\star}} | \\
	& = \left\{ \sum_{j\geq 0} \left[(2^{dj})^{\frac{\gamma}{d}+\frac{1}{2}-\frac{1}{p}} \| f_{j\cdot} \|_{p} \right]^q \right\}^{1/q} \left\{ \sum_{j\geq 0} \left[(2^{-dj})^{\frac{\gamma}{d}+\frac{1}{2}-\frac{1}{p}} \| u_{j\cdot} - v_{j\cdot} \|_{p_{\star}} \right]^{q_{\star}} \right\}^{1/q_{\star}} \enspace.
\end{align*}
Take $p=q=\infty$ (in this case $p_{\star} = q_{\star} = 1$), we know that the IPM can be regarded as a type of Besov space semi-norm
\begin{align*}
	d_{\cF_\gamma}(\mu, \nu) = \sum_{j\geq 0} (2^{-dj})^{\frac{\gamma}{d}+\frac{1}{2}} \sum_{k=0}^{2^j-1} |u_{jk} - v_{jk}|.
\end{align*}

\paragraph{Step 2: composite hypothesis and prior construction.}
Next we are going to construct two priors on $\nu$ such that the difference
\begin{align*}
	|\E_{\nu \sim \cP_0} d_{\cF_\gamma}(\mu, \nu) - \E_{\nu \sim \cP_1} d_{\cF_\gamma}(\mu, \nu) |
\end{align*}
is large, while at the same time one can not distinguish the following two distributions 
\begin{align}
	p_0(Y_1, \ldots Y_n) =  \E_{\nu \sim \cP_0} \left[ \prod_{i=1}^n \rho_\nu(Y_i) \right],~ p_1(Y_1, \ldots Y_n) =  \E_{\nu \sim \cP_1} \left[ \prod_{i=1}^n \rho_\nu(Y_i) \right] .
\end{align}
Here $\cP_0, \cP_1$ are two prior distributions on $\nu$ which we will construct. Consider $\mu$ to be the same distribution under the null $H_0$ and the alternative $H_1$.
Set two values $K, \tau$ to be used in the construction
\begin{align}
	K \asymp \frac{\log n}{\log \log n}, ~\tau \asymp 1.
\end{align}
The choice will be apparent in the latter part of the proof. The following prior construction is inspired by \cite{lepski1999estimation}, where they study the estimation of functionals under the Gaussian white noise model. This prior was also used in \cite{cai2011testing} for studying non-smooth functional estimation in Gaussian sequence models.
\begin{proposition}[\cite{lepski1999estimation}, Proposition 4.2]
	\label{prop:hardest-priors}
	For any given positive integer $K$ and $\tau \in \mathbb{R}_{\geq 0}$, there exist two symmetric probability measures $q_0$ and $q_1$ on $[-\tau, \tau]$ such that 
	 \begin{align}
	 	\int_{-\tau}^\tau t^l q_0(dt) = \int_{-\tau}^\tau t^l q_1(dt), ~~ l = 0, 1, \ldots, 2K, \\
		\int_{-\tau}^\tau |t| q_1(dt) - \int_{-\tau}^\tau |t| q_0(dt) = 2\kappa \cdot K^{-1} \tau,
	 \end{align}
	 where $\kappa$ is some constant depending on $K$ only.
\end{proposition}

Now let's construct $\cP_0$ and $\cP_1$ as follows. Take $\mu \sim {\rm Unif}([0,1]^d)$. Choose $J \in \mathbb{N}$ such that $2^{dJ} \asymp n^{\frac{1}{1+2\beta/d}}$, first we embed a parametrized class of densities into $\cG_\beta$
\begin{align}
	\label{eq:nu_theta}
	\frac{d\nu_{\theta}}{dx} := \mu(x) + \frac{1}{\sqrt{n}} \sum_{k=0}^{2^{dJ}-1} \theta_k h_{Jk}(x)
\end{align}
with each $\theta_{k} \in [-\tau, \tau]$ for all $k$.
Now we show that the construction lies inside the space of interest, i.e., $\nu_\theta \in \cG_{\beta}$.
First observe that for the wavelet basis that satisfy the regularity condition $\int_{\Omega} h_{jk} d\mu = 0$, we have $\int_{\Omega} d \nu_\theta  = 1$ and $\| d \nu_\theta/dx \|_{L_\infty} \geq 1 - \sqrt{2^{dJ}/n} >0$. Hence $\nu_\theta$ is a valid probability measure.
Let's then verify that the density $\rho_{\nu_\theta} \in \sB^{\beta, \infty}_{\infty}$ lies in the Besov space, since
\begin{align}
	\frac{1}{\sqrt{n}} |\theta_k| \leq (2^{dJ})^{-(\frac{\beta}{d}+\frac{1}{2})}, ~~\forall k.
\end{align} 
For any $\gamma \geq 0$, it is then easy to verify via Step 1 that
\begin{align}
	d_{\cF_\gamma}(\mu, \nu_{\theta})&:=  (2^{-dJ})^{\frac{\gamma}{d}+\frac{1}{2}} \frac{1}{\sqrt{n}}\sum_{k=0}^{2^{dJ}-1} |\theta_k| \nonumber\\
	&=(2^{-dJ})^{\frac{\gamma}{d}+\frac{1}{2}} (2^{dJ})^{-(\frac{\beta}{d}+\frac{1}{2})} \sum_{k=0}^{2^{dJ}-1} |\theta_k| \nonumber\\
	&=(2^{-dJ})^{\frac{\beta+\gamma}{d}} \frac{1}{2^{dJ}}\sum_{k=0}^{2^{dJ}-1} |\theta_k|.
\end{align}

Making use of the probability measures $q_0$ and $q_1$ on $[-\tau, \tau]$ claimed by Proposition~\ref{prop:hardest-priors}, we define a collection of measures $$\cS_{0} := \{\nu_\theta: \theta_k \sim q_0~ i.i.d.~ \text{for each}~ k\in [2^{dJ}] \}.$$ Then $\cP_0$ can be viewed as an uniform prior over this set $\cS_{0}$. Similar construction holds for $\cP_1$ via $q_1$.

\paragraph{Step 3: polynomials and matching moments.}

Remark that due to the separation of support for wavelets (localized property), i.e., $h_{Jk}(x) h_{Jk'}(x) = 0$ for $k \neq k'$, we have the equivalent expression as in \eqref{eq:nu_theta}
\begin{align}
	\label{eq:mu_theta_prod}
	\frac{d\nu_{\theta}}{dx} = \prod_{k=0}^{2^{dJ}-1} (1+ \theta_k n^{-1/2} h_{Jk}(x)) \enspace.
\end{align}
Use $\theta\sim q_0^{\otimes 2^{dJ}}$ to denote that $\theta_k \sim q_0$ i.i.d. for all $k \in [2^{dJ}]$, we know
\begin{align}
	\label{eq:sum-prod}
	p_0(Y_1,\ldots,Y_n) &= \E_{\theta \sim q_0^{\otimes 2^{dJ}}}  \prod_{i=1}^n \rho_\theta(Y_i) = \E_{\theta \sim q_0^{\otimes 2^{dJ}}}  \prod_{i=1}^n \prod_{k=0}^{2^{dJ}-1} (1+\theta_k n^{-1/2} h_{Jk}(Y_i)) \nonumber \\
	&= \E_{\theta \sim q_0^{\otimes 2^{dJ}}}  \prod_{k=0}^{2^{dJ}-1} \prod_{i=1}^n (1+\theta_k n^{-1/2} h_{Jk}(Y_i)) \nonumber \quad \text{by \eqref{eq:mu_theta_prod}}\\
	&= \prod_{k=1}^{2^{dJ}} \E_{\theta_k \sim q_0} \prod_{i=1}^n (1+\theta_k n^{-1/2} h_{Jk}(Y_i)) \enspace.
\end{align}
Remark that we can not further interchange the ordering of $\E_{\theta_k}$ and $\prod_{i=1}^n$, since the mixture is on data distributions $(Y_1, \ldots, Y_n)$ jointly.

Let's introduce the polynomial $f(\theta_k;h_{jk}(Y_1),\ldots,h_{jk}(Y_n))$ in $\theta_k$ (and $h_{Jk}(Y_i)$) with degree at most $n$ appearing in the above expression, which will be used extensively in the next step,
\begin{align}
	&f(\theta_k;h_{jk}(Y_1),\ldots,h_{jk}(Y_n)):=  \prod_{i=1}^n (1+\theta_k \frac{h_{Jk}(Y_i)}{\sqrt{n}}) \label{eq:polyn}\\
	&=\sum_{l=0}^{n} \theta_k^l \frac{\sum_{i_1<\ldots<i_l} h_{Jk}(Y_{i_1})\ldots h_{Jk}(Y_{i_l})}{n^{l/2}} =: \sum_{l=0}^{n} \theta_k^l \frac{H^{(l)}_{Jk}(Y_1, \ldots, Y_n)}{n^{l/2}} \nonumber \enspace.
\end{align}
Here $H^{(l)}_{JK}(Y_1, \ldots, Y_n)$ is a sum of monomials of order $l$, i.e., $\binom{n}{l}$ terms with each of the form $h_{Jk}(Y_{i_1})\ldots h_{Jk}(Y_{i_l})$.
Denote $f^{[\leq K]}, f^{[> K]}$ to be the corresponding truncated polynomial according to the degree. In this convenient notation, we know
\begin{align}
	p_0(Y_1,\ldots,Y_n) = \prod_{k \in [2^{dJ}]} \E_{\theta_k \sim q_0}f(\theta_k;h_{Jk}(Y_1),\ldots,h_{Jk}(Y_n)) .
\end{align}
Later, we shall use the following properties of the polynomial $f$ of degree at most $n$,
\begin{align}
	\forall \theta_k, ~~\int_{\cY^{\otimes n}} f(\theta_k;h_{Jk}(y_1),\ldots,h_{Jk}(y_n)) dy_1\ldots dy_n= 1 \enspace.
\end{align}
And the following property according to $q_0$ and $q_1$ constructed in Proposition~\ref{prop:hardest-priors}: $\forall y_1,\ldots,y_n$
\begin{align*}
	&\E_{\theta_k \sim q_1} f(\theta_k;h_{Jk}(y_1),\ldots,h_{Jk}(y_n)) - \E_{\theta_k \sim q_0} f(\theta_k;h_{Jk}(y_1),\ldots,h_{Jk}(y_n)) \\
	&= \int_{[-\tau,\tau]} f^{[>2K]}(\theta_k;h_{Jk}(y_1),\ldots,h_{Jk}(y_n)) (q_1-q_0)(d\theta_k) \enspace.
\end{align*}

\paragraph{Step 4: total variation, telescoping and the sum-product trick.}
When there is no confusion, we use $f(\theta_k;h_{Jk}(y^{\otimes n}))$ to abbreviate $f(\theta_k; h_{Jk}(y_1),\ldots, h_{Jk}(y_n))$. Recall \eqref{eq:sum-prod}, we have
\begin{align*}
	{\rm TV}(p_1, p_0)&:= \frac{1}{2} \int_{\Omega^{\otimes n}} \left|p_1(y_1, \ldots, y_n) - p_0(y_1, \ldots, y_n) \right| dy_1\ldots dy_n  \\
	&= \frac{1}{2} \int_{\Omega^{\otimes n}}  \left|   \prod_{k\in [2^{dJ}]} \E_{\theta_k \sim q_1} f(\theta_k;h_{Jk}(y^{\otimes n})) - \prod_{k \in [2^{dJ}]} \E_{\theta_k \sim q_0}f(\theta_k;h_{Jk}(y^{\otimes n})) \right|  dy_1\ldots dy_n  .\\
\end{align*}

We claim that the following telescoping Lemma holds. The proof can be seen clearly through writing the left hand side as a telescoping sum and using the triangle inequality.
\begin{proposition}[Telescoping]
	For $N \in \mathbb{N}, N\geq 2$, and $a_i, b_i \geq 0, 1 \leq i \leq N$,
	\begin{align}
		|\prod_{k\in[1, N]} a_k - \prod_{k \in [1, N]} b_k| \leq \sum_{i\in [1, N]} |a_i - b_i| \cdot \prod_{k\in [1, i)} b_k  \cdot \prod_{k\in (i, N]} a_k \enspace.
	\end{align}
\end{proposition}

To make use of the above Lemma, define 
\begin{align}
	a_k(h_{Jk}(y_1),\ldots, h_{Jk}(y_n)) := \E_{\theta_k \sim q_1} f(\theta_k;h_{Jk}(y^{\otimes n})) \\
	b_k(h_{Jk}(y_1),\ldots, h_{Jk}(y_n)) := \E_{\theta_k \sim q_0} f(\theta_k;h_{Jk}(y^{\otimes n}))
\end{align}
Using the the above telescoping proposition, we claim
\begin{align}
	{\rm TV}(p_1, p_0)&\leq \sum_{k \in [2^{dJ}]} \int |a_k - b_k| \cdot  \left( \prod_{k'\in [0, k)} b_{k'}  \prod_{k''\in (k, 2^{dJ}-1]} a_{k''}  dy_1\ldots dy_n   \right) \nonumber \\
	&= \sum_{k \in [2^{dJ}]} \E_{\substack{\theta_{k'} \sim q_0, k' \in [0, k)\\
	\theta_{k''} \sim q_1, k'' \in (k, 2^{dJ}-1]} } \E_{Y_1,\ldots,Y_n \sim \nu_{\theta_{-k}}} |a_k(h_{Jk}(Y_1),\ldots, h_{Jk}(Y_n)) - b_k(h_{Jk}(Y_1),\ldots, h_{Jk}(Y_n))| . \label{eq:cruc}
\end{align}
The reasoning behind the last line is as follows.
Firstly, we need to define a tilted measure $\nu_{\theta_{-k}}$ without the influence of the $k$-the coordinate $\theta_k$,
\begin{align}
	\rho_{\nu_{\theta_{-k}}}(x) = \frac{d \nu_{\theta_{-k}}}{dx} := 1+ \frac{1}{\sqrt{n}} \sum_{\substack{k'\neq k\\ 0 \leq k' \leq 2^{dJ}-1 }} \theta_{k'} h_{J k'}(x) =\prod_{\substack{k'\neq k\\ 0 \leq k' \leq 2^{dJ}-1 }} \left( 1+ \frac{1}{\sqrt{n}} \theta_{k'} h_{J k'}(x) \right) \enspace.
\end{align}
From the properties established in Step 3, one can verify that
\begin{align}
	\E_{\substack{\theta_{k'} \sim q_0, k' \in [0, k)\\ \theta_{k''} \sim q_1, k'' \in (k, 2^{dJ}-1]} } \prod_{i=1}^n \rho_{\nu_{\theta_{-k}}}(y_i) & = \E_{\substack{\theta_{k'} \sim q_0, k' \in [0, k)\\\theta_{k''} \sim q_1, k'' \in (k, 2^{dJ}-1]} } \prod_{i=1}^n \prod_{\substack{k'\neq k\\ 0 \leq k' \leq 2^{dJ}-1 }} \left( 1+ \frac{1}{\sqrt{n}} \theta_{k'} h_{J k'}(y_i) \right) \nonumber \\
		& = \E_{\substack{\theta_{k'} \sim q_0, k' \in [0, k)\\\theta_{k''} \sim q_1, k'' \in (k, 2^{dJ}-1]} }  \prod_{\substack{k'\neq k\\ 0 \leq k' \leq 2^{dJ}-1 }} \prod_{i=1}^n \left( 1+ \frac{1}{\sqrt{n}} \theta_{k'} h_{J k'}(y_i) \right) \nonumber  \\
		&=  \prod_{k'\in [0, k)} b_{k'}  \prod_{k''\in (k, 2^{dJ}-1]} a_{k''} \enspace.
\end{align}
Now we have proved \eqref{eq:cruc}, since by using Fubini's theorem,
\begin{align*}
	&\int |a_k - b_k| \cdot  \left( \prod_{k'\in [0, k)} b_{k'}  \prod_{k''\in (k, 2^{dJ}-1]} a_{k''}  dy_1\ldots dy_n   \right) \\
	&= \E_{\substack{\theta_{k'} \sim q_0, k' \in [0, k)\\ \theta_{k''} \sim q_1, k'' \in (k, 2^{dJ}-1]} }  \int |a_k(h_{Jk}(y_1),\ldots, h_{Jk}(y_n)) - b_k(h_{Jk}(y_1),\ldots, h_{Jk}(y_n))| \prod_{i=1}^n \rho_{\nu_{\theta_{-k}}}(y_i) dy_1\ldots dy_n \enspace.
\end{align*}

Let's analyze the term
\begin{align*}
	\E_{Y_1,\ldots,Y_n \sim \nu_{\theta_{-k}}} |a_k(h_{Jk}(Y_1),\ldots, h_{Jk}(Y_n)) - b_k(h_{Jk}(Y_1),\ldots, h_{Jk}(Y_n))| 
\end{align*}
where $Y_1, \ldots Y_n$ are i.i.d. sampled from a measure $\nu_{\theta_{-k}}$.
We emphasize that $\nu_{\theta_{-k}}$ agrees with the uniform measure $\mu$ on the domain associated with $h_{Jk}(x)$.
Due to the separation of support for the wavelet basis, we know that the random variables
\begin{align*}
	h_{J k}(Y_i)
\end{align*}
are only determined by $\nu_{\theta_{-k}}$ restricted to the domain of $h_{Jk}$. Equivalently, the distributions of $h_{J k}(Y)$'s are the same when $Y\sim \nu_{\theta_{-k}}$ and $Y \sim \mu$.
Hence for $Y_1, \ldots, Y_n \sim \nu_{\theta_{-k}}$, 
\begin{align*}
	&\E_{Y_1,\ldots,Y_n \sim \nu_{\theta_{-k}}} |a_k(h_{Jk}(Y_1),\ldots, h_{Jk}(Y_n)) - b_k(h_{Jk}(Y_1),\ldots, h_{Jk}(Y_n))|\\
	& = \E_{Y_1,\ldots,Y_n \sim \mu} |a_k(h_{Jk}(Y_1),\ldots, h_{Jk}(Y_n)) - b_k(h_{Jk}(Y_1),\ldots, h_{Jk}(Y_n))| \enspace.
\end{align*}
Now one can directly bound the TV metric between the complex sum-product distribution $p_0$ and $p_1$ defined in \eqref{eq:sum-prod},
\begin{align}
	\label{eq:teles-later}
	2{\rm TV}(p_1, p_0)& \leq \sum_{k=0}^{2^{dJ}-1} \E_{Y_1,\ldots,Y_n \sim \mu} |a_k(h_{Jk}(Y_1),\ldots, h_{Jk}(Y_n)) - b_k(h_{Jk}(Y_1),\ldots, h_{Jk}(Y_n))|  \nonumber \\
	& = \sum_{k=0}^{2^{dJ}-1} \int \left| \E_{\theta_k \sim q_1} f(\theta_k;h_{Jk}(y^{\otimes n})) - \E_{\theta_k \sim q_0}f(\theta_k;h_{Jk}(y^{\otimes n}))\right| dy_1 \ldots dy_n .
\end{align}

\paragraph{Step 5: $\ell_2$ bound.}
In this section, we are going to bound, for a fixed $k$, the following expression using the properties of the $q_1$ and $q_0$ constructed with matching moments up to $2K$ (claimed by Proposition~\ref{prop:hardest-priors}),
\begin{align*}
	\int \left| \E_{\theta_k \sim q_1} f(\theta_k;h_{Jk}(y^{\otimes n})) - \E_{\theta_k \sim q_0}f(\theta_k;h_{Jk}(y^{\otimes n}))\right| dy_1 \ldots dy_n \enspace.
\end{align*}
First, observe the $\ell_2$ bound
\begin{align}
	\int |g_1 -g_2| d\mu \leq \left( \int (g_1 - g_2)^2 d\mu \right)^{1/2} \enspace.
\end{align}
Let's bound the $\ell_2$ form, which takes the form
\begin{align}
	&\int \left( \E_{\theta_k \sim q_1} f(\theta_k;h_{Jk}(y^{\otimes n})) - \E_{\theta_k \sim q_0}f(\theta_k;h_{Jk}(y^{\otimes n}))\right)^2 dy_1 \ldots dy_n \label{eq:l2}\\
	&= \E_{\theta, \theta'\sim q_1} \int f(\theta;h_{Jk}(y^{\otimes n})) f(\theta';h_{Jk}(y^{\otimes n})) d y^{\otimes n} + \E_{\omega, \omega'\sim q_0} \int f(\omega;h_{Jk}(y^{\otimes n})) f(\omega';h_{Jk}(y^{\otimes n})) d y^{\otimes n} \nonumber\\
	& \quad - 2 \E_{\theta \sim q_1, \omega\sim q_0}\int f(\theta;h_{Jk}(y^{\otimes n})) f(\omega;h_{Jk}(y^{\otimes n})) d y^{\otimes n} \nonumber \enspace.
\end{align}
Note now each $f(\theta;h_{Jk}(y^{\otimes n}))f(\theta';h_{Jk}(y^{\otimes n}))$ for fixed $\theta, \theta'$ takes the following product form
\begin{align*}
	f(\theta;h_{Jk}(y^{\otimes n}))f(\theta';h_{Jk}(y^{\otimes n})) = \prod_{i=1}^n \left(1+ (\theta + \theta') \frac{h_{Jk}(Y_i)}{\sqrt{n}} + \theta \theta' \frac{h^2_{Jk}(Y_i)}{n}  \right)
\end{align*}
and
\begin{align*}
	\int f(\theta;h_{Jk}(y^{\otimes n})) f(\theta';h_{Jk}(y^{\otimes n})) d y^{\otimes n} &= \left( 1 + \theta\theta' \frac{\int h^2_{Jk}(y) dy}{n} \right)^n \\
	&= \left( 1 + \theta\theta' \frac{1}{n} \right)^n \enspace.
\end{align*}

Therefore we have for \eqref{eq:l2} 
\begin{align*}
	\eqref{eq:l2}&= \E_{\theta, \theta'\sim q_1} \left[\left( 1 + \theta\theta' \frac{1}{n} \right)^n \right] + \E_{\omega, \omega'\sim q_0} \left[\left( 1 + \omega\omega' \frac{1}{n} \right)^n \right] -2 \E_{\theta \sim q_1, \omega\sim q_0} \left[\left( 1 + \theta\omega \frac{1}{n} \right)^n \right] \\
	&= \sum_{l=1}^{\lfloor n/2\rfloor}  \left(\E_{\theta, \theta'\sim q_1} [(\theta \theta')^{2l}] + \E_{\omega, \omega'\sim q_0} [(\omega \omega')^{2l}] - 2 \E_{\theta \sim q_1, \omega\sim q_0} [(\theta \omega)^{2l}] \right) \frac{\binom{n}{2l}}{n^{2l}} \\
	&= \sum_{l=1}^{\lfloor n/2\rfloor}  \left( \left(\E_{q_1} [\theta^{2l}] \right)^2 + \left(\E_{q_0} [\theta^{2l}] \right)^2 - 2 \E_{q_1} [\theta^{2l}]  \E_{q_0} [\theta^{2l}] \right) \frac{\binom{n}{2l}}{n^{2l}} 
\end{align*}
Recall the crucial property that for all $l\leq K$, we know
\begin{align}
	\E_{\theta \sim q_1} [\theta^{2l}] = \E_{\theta \sim q_0} [\theta^{2l}] ~~ \Rightarrow~~
	\left(\E_{q_1} [\theta^{2l}] \right)^2 + \left(\E_{q_0} [\theta^{2l}] \right)^2 - 2 \E_{q_1} [\theta^{2l}]  \E_{q_0} [\theta^{2l}] = 0
\end{align}
therefore the above summation equals
\begin{align*}
	\eqref{eq:l2}&=\sum_{l=K+1}^{\lfloor n/2\rfloor} \left( \left(\E_{q_1} [\theta^{2l}] \right)^2 + \left(\E_{q_0} [\theta^{2l}] \right)^2 - 2 \E_{q_1} [\theta^{2l}]  \E_{q_0} [\theta^{2l}] \right) \frac{\binom{n}{2l}}{n^{2l}} \\
	&\leq \sum_{l=K+1}^{\lfloor n/2\rfloor}  4\tau^{4l} \frac{1}{(2l)!} \\
	& \precsim 4\frac{\tau^{4K}}{(2K)!} \exp(\tau^4) \enspace.
\end{align*}

Assemble the two bounds, we have
\begin{align}
	\label{eq:crucial-ineq}
	&\int \left| \E_{\theta_k \sim q_1} f(\theta_k;h_{Jk}(y^{\otimes n})) - \E_{\theta_k \sim q_0}f(\theta_k;h_{Jk}(y^{\otimes n}))\right| dy_1 \ldots dy_n \\
	&\leq 2 \frac{\tau^{2K}}{\sqrt{(2K)!}} \exp(\tau^4/2)
\end{align}

\paragraph{Step 6: combine all pieces.}

Now continuing \eqref{eq:teles-later}, we have
\begin{align*}
	2{\rm TV}(p_1, p_0)& \leq \sum_{k=0}^{2^{dJ}-1} \E_{Y_1,\ldots,Y_n \sim \mu} |a_k(h_{Jk}(Y_1),\ldots, h_{Jk}(Y_n)) - b_k(h_{Jk}(Y_1),\ldots, h_{Jk}(Y_n))| \\
	& = \sum_{k=0}^{2^{dJ}-1} \int \left| \E_{\theta_k \sim q_1} f(\theta_k;h_{Jk}(y^{\otimes n})) - \E_{\theta_k \sim q_0}f(\theta_k;h_{Jk}(y^{\otimes n}))\right| dy_1 \ldots dy_n \\
	&\leq 2^{dJ}  \cdot 2 \frac{\tau^{2K}}{\sqrt{2K}!} \exp(\tau^4/2) \precsim \exp (c \log n - K \log K) \enspace.
\end{align*}
Therefore by taking $K = \frac{c}{2} \frac{\log n}{\log \log n}$, we know
\begin{align}
	2{\rm TV}(p_1, p_0) \leq \exp(-\frac{c}{2}\log n) \leq n^{-c/2}. \label{eq:tv-bd}
\end{align}

By the construction of the composite hypothesis, we have
\begin{align*}
	& \quad |\E_{\nu_{\theta} \sim \cP_0} d_{\cF_\gamma}(\mu, \nu_{\theta}) - \E_{\nu_{\theta} \sim \cP_1} d_{\cF_\gamma}(\mu, \nu_{\theta})| \\
	&=  (2^{-dJ})^{-\frac{\beta+\gamma}{d}} \cdot \left| \E_{\nu_{\theta} \sim \cP_0} \left[ \frac{1}{2^{dJ}} \sum_{k\in [2^{dJ}]}|\theta_k| \right] - \E_{\nu_{\theta} \sim \cP_1} \left[ \frac{1}{2^{dJ}} \sum_{k\in [2^{dJ}]}|\theta_k| \right]   \right|  \\
	& = n^{-\frac{\beta+\gamma}{2\beta+d}} \cdot \left| \E_{\theta \sim q_0} [|\theta|] -\E_{\theta \sim q_1}[ |\theta|] \right| \\
	& \geq n^{-\frac{\beta+\gamma}{2\beta+d}} \cdot 2\kappa K^{-1} \tau  \asymp n^{-\frac{\beta+\gamma}{2\beta+d}} \cdot \frac{\log \log (n)}{\log(n)} \enspace.
\end{align*}

Denote $\cD_n$ to be the collection of data $(Y_1, \cdots, Y_n)$, which is drawn from the distribution $Pr(y^{\otimes n} |\theta) := \prod_{i=1}^n \rho_{\nu_\theta}(y_i)$.
For any functional of $\theta$, and for any estimator based on $n$-i.i.d. samples, we know
\begin{align*}
	\sup_{\nu_\theta}\E_{\cD_n \sim Pr(y^{\otimes n} |\theta)} |\hat{T}_n - F(\theta)| &\geq \E_{\theta \sim Q_0} \E |\hat{T}_n - F(\theta)| \\
	&\geq \E_{\theta \sim Q_0} \E_{\cD_n \sim Pr(y^{\otimes n} |\theta)} |\hat{T}_n - \E_{\theta\sim Q_0}F(\theta)| - \delta_{Q_0}
\end{align*}
where $\delta_{Q_0}:=\E_{\theta \sim Q_0} |\E_{\theta\sim Q_0}F(\theta) - F(\theta)| $. Here $Q_0$ is some prior distribution on $\theta$.
Repeat the same argument for $Q_1$,
and by Le Cam's argument on two composite hypothesis
\begin{align*}
	\sup_{\nu_\theta} \E |\hat{T}_n - F(\theta)|& \geq \frac{1}{2}\left( \E_{\theta \sim Q_0} \E_{\cD_n \sim Pr(y^{\otimes n} |\theta)} |\hat{T}_n - \E_{\theta\sim Q_0}F(\theta)| + \E_{\theta \sim Q_1} \E_{\cD_n \sim Pr(y^{\otimes n} |\theta)} |\hat{T}_n - \E_{\theta\sim Q_1}F(\theta)| \right) - \frac{\delta_{Q_0} + \delta_{Q_1}}{2} \\
	& = \frac{1}{2}\left( \E_{\cD_n \sim p_0}|\hat{T}_n - \E_{\theta\sim Q_0}F(\theta)| + \E_{\cD_n \sim p_1} |\hat{T}_n - \E_{\theta\sim Q_1}F(\theta)| \right) - \frac{\delta_{Q_0} + \delta_{Q_1}}{2} \\
	& \geq \frac{|\E_{\theta\sim Q_0}F(\theta) - \E_{\theta\sim Q_1}F(\theta)|}{4} \left( P_{0}(T = 1) + P_{1}(T = 0)\right) - \frac{\delta_{Q_0} + \delta_{Q_1}}{2} \\
	& \geq \frac{|\E_{\theta\sim Q_0}F(\theta) - \E_{\theta\sim Q_1}F(\theta)|}{4} \int p_0(y^{\otimes n}) \wedge p_1(y^{\otimes n}) dy^{\otimes n} -  \frac{\delta_{Q_0} + \delta_{Q_1}}{2} \\
	& = \frac{|\E_{\theta\sim Q_0}F(\theta) - \E_{\theta\sim Q_1}F(\theta)|}{4} (1 - d_{TV}(p_0, p_1)) -    \frac{\delta_{Q_0} + \delta_{Q_1}}{2} 
\end{align*}
where the posterior distribution $p_i(y^{\otimes n}) = \int Pr(y^{\otimes n} |\theta) Q_i(d\theta)$, for $i=0,1$. Here the test $T=1$ if and only if $\hat{T}_n$ is closer to $\E_{\theta\sim Q_1}F(\theta)$.
In our case, 
$$
F(\theta) := d_{\cF_\gamma}(\mu, \nu_{\theta}) = (2^{-dJ})^{-\frac{\beta+\gamma}{d}} \left[ \frac{1}{2^{dJ}} \sum_{k\in [2^{dJ}]}|\theta_k| \right] \enspace,
$$
hence we know
\begin{align*}
	|\E_{\theta\sim Q_0}F(\theta) - \E_{\theta\sim Q_1}F(\theta)| &= |\E_{\nu_{\theta} \sim \cP_0} d_{\cF_\gamma}(\mu, \nu_{\theta}) - \E_{\nu_{\theta} \sim \cP_1} d_{\cF_\gamma}(\mu, \nu_{\theta})| \\
	&\succsim n^{-\frac{\beta+\gamma}{2\beta+d}} \cdot \frac{\log \log (n)}{\log(n)} \\
	 1 - d_{TV}(p_0, p_1) &\geq 1 -  n^{-c/2} \quad \text{by \eqref{eq:tv-bd}}\\
	 \frac{\delta_{Q_0} + \delta_{Q_1}}{2} &\precsim n^{-\frac{\beta+\gamma}{2\beta+d}} \frac{1}{\sqrt{2^{dJ}}} \ll n^{-\frac{\beta+\gamma}{2\beta+d}} \cdot \frac{\log \log (n)}{\log(n)} \enspace.
\end{align*}
Therefore we have
\begin{align}
		\inf_{\widehat{T}_n} \sup_{\nu \in \sC^{\beta}} \E|\widehat{T}_n - d_{\cF_\gamma}(\mu, \nu)| \succsim  n^{-\frac{\beta+\gamma}{2\beta+d}} \cdot \frac{\log \log (n)}{\log(n)} \enspace.
\end{align}


\subsection{Proof of the Upper Bound}

The upper bound can be obtained through similar derivations as in \cite{liang2018well, singh2018nonparametric, weed2019estimation}. We include here for completeness. 

The estimator is of the plug-in form, with 
\begin{align}
	d_{\cF_{\gamma}}(\widetilde \mu_m, \widetilde \nu_n) := \sup_{f \in \cF_{\gamma}} |\int f d \widetilde{\mu}_m - \int f d\widetilde{\nu}_n |
\end{align}
where $\widetilde \mu_m$, and $\widetilde \nu_n$ are smoothed empirical measures based on truncation on Wavelets.
It is clear that
\begin{align}
	|d_{\cF_{\gamma}}(\widetilde \mu_m, \widetilde \nu_n) - d_{\cF_{\gamma}}(\mu, \nu)| \leq \sup_{f \in \cF_{\gamma}} |\int f d \widetilde{\mu}_m - \int f d\mu | + \sup_{f \in \cF_{\gamma}} |\int f d \widetilde{\nu}_n - \int f d\nu |.
\end{align}

Now let's bound $\sup_{f \in \cF_{\gamma}} |\int f d \widetilde{\nu}_n - \int f d\nu |$ via expanding under the Wavelet basis. Denote $\widehat{\E}[h_{jk}]:=1/n\sum_{i=1}^n h_{jk}(Y_i)$,  the smoothed empirical estimate $\widetilde\nu_n$ is defined as
\begin{align}
	\frac{d \widetilde\nu_n}{dx} := \sum_{j=0}^{J} \sum_{k=0}^{2^{dj}-1} \widehat{\E}[h_{jk}] h_{jk}(x) \enspace.
\end{align}
Expand $f(x) = \sum_{j\geq 0} \sum_{k=0}^{2^{dj}-1} f_{jk} h_{jk}(x)$, we have
\begin{align*}
	 & \sup_{f \in \cF_{\gamma}} |\int f d \widetilde{\nu}_n - \int f d\nu | \leq \sup_{f \in  \sB^{\gamma,\infty}_{\infty}(1)} |\int f d \widetilde{\nu}_n - \int f d\nu | \\
	 & = \sup_{f \in  \sB^{\gamma,\infty}_{\infty}(1)} | \sum_{j \geq 0}^{J}  \sum_{k=0}^{2^{dj}-1} f_{jk} (\widehat{\E}[h_{jk}] - \E[h_{jk}])| +  \sup_{f \in  \sB^{\gamma,\infty}_{\infty}(1)} | \sum_{j >J}  \sum_{k=0}^{2^{dj}-1} f_{jk}  \E[h_{jk}] |
\end{align*}
For the first term, since $f \in  \sB^{\gamma,\infty}_{\infty}(1) \Rightarrow \forall j,k, ~|f_{jk}| \leq (2^{-dj})^{\frac{\gamma}{d}+\frac{1}{2}}$
\begin{align*}
	&\E \sup_{f \in  \sB^{\gamma,\infty}_{\infty}(1)} | \sum_{j \geq 0}^{J}  \sum_{k=0}^{2^{dj}-1} f_{jk} (\widehat{\E}[h_{jk}] - \E[h_{jk}])| \leq \sum_{j\geq 0}^J (2^{-dj})^{\frac{\gamma}{d}+\frac{1}{2}} \sum_{k=0}^{2^{dj}-1} \E|\widehat{\E}[h_{jk}] - \E[h_{jk}]|\\
	&\leq \sum_{j\geq 0}^J (2^{-dj})^{\frac{\gamma}{d}+\frac{1}{2}} \sum_{k=0}^{2^{dj}-1} (\E|\widehat{\E}[h_{jk}] - \E[h_{jk}]|^2)^{1/2} \quad \text{since $\sqrt{\E [Z]} \geq \E [\sqrt{Z}]$ for $Z \geq 0$}\\
	&\precsim \sum_{j\geq 0}^J (2^{-dj})^{\frac{\gamma}{d}+\frac{1}{2}} 2^{dj} \frac{1}{\sqrt{n}} \asymp \frac{1}{\sqrt{n}}(2^{dJ})^{\frac{1}{2} - \frac{\gamma}{d}} 
\end{align*}
for $d\geq 2\gamma$.

For the second term, recall $\E_{Y\sim \nu} [h_{jk}(Y)] = \langle d\nu/dx, h_{jk} \rangle =: v_{jk}$. Due to the fact that
\begin{align}
	d \nu/dx \in \sB^{\beta, \infty}_{\infty}  \Rightarrow \forall j,k, ~|v_{jk}| \leq (2^{-dj})^{\frac{\beta}{d}+\frac{1}{2}}\\
	f \in \sB^{\gamma, \infty}_{\infty} \Rightarrow \forall j,k, ~|f_{jk}| \leq (2^{-dj})^{\frac{\gamma}{d}+\frac{1}{2}}
\end{align}
\begin{align*}
	& \E \sup_{f \in \sB^{\gamma,\infty}_{\infty}} | \sum_{j >J}  \sum_{k=0}^{2^{dj}-1} f_{jk} \E[h_{jk}] | = \E \sup_{f \in \sB^{\gamma,\infty}_{\infty}} | \sum_{j >J}  \sum_{k=0}^{2^{dj}-1} f_{jk} v_{jk} | \\
	& \leq  \sum_{j >J}  \sum_{k=0}^{2^{dj}-1} (2^{-dj})^{\frac{\gamma}{d}+\frac{1}{2}} (2^{-dj})^{\frac{\beta}{d}+\frac{1}{2}} \\
	& \leq (2^{dJ})^{-\frac{\beta+\gamma}{d}} \enspace .
\end{align*}

Balancing the two terms, we have
\begin{align}
	& \sup_{\nu \in \cG_\beta} \E \sup_{f \in \cF_{\gamma}} |\int f d \widetilde{\nu}_n - \int f d\nu | \precsim \frac{1}{\sqrt{n}}(2^{dJ})^{\frac{1}{2} - \frac{\gamma}{d}}  + (2^{dJ})^{-\frac{\beta+\gamma}{d}}\\
	& \asymp n^{-\frac{\beta+\gamma}{2\beta+d}}, ~~\text{with $2^{dJ} \asymp n^{\frac{1}{2\beta/d+1}}$}\enspace.
\end{align}
Put everything together, we know
\begin{align}
	\E|d_{\cF_{\gamma}}(\widetilde \mu_m, \widetilde \nu_n) - d_{\cF_{\gamma}}(\mu, \nu)| \leq (n \wedge m)^{-\frac{\beta+\gamma}{2\beta+d}}.
\end{align}

\bibliographystyle{chicago}
\bibliography{ref}

\end{document}